\documentclass{article}
\usepackage{amsfonts,amsmath,amsthm}
\theoremstyle{plain}
\newtheorem{theo}{Theorem}
\newtheorem{cor}[theo]{Corollary}
\newtheorem{prop}[theo]{Proposition}
\theoremstyle{definition}
\newtheorem*{defi}{Definition}
\newtheorem*{rem}{Remark}
\def\proof{\par\medskip\noindent{\em Proof}:\ }
\begin{document}
\title{A distinguished geometry perspective on\\ multi-time affine quadratic Lagrangians}
\author{Mircea Neagu}
\date{}
\maketitle
\begin{abstract}
For a space endowed with a general quadratic multi-time Lagrangian and an associated non-linear connection,
    the paper constructs the main Riemann-Lagrange distinguished geometric objects (linear connection, torsion
    and curvature).
\end{abstract}


\section{Introduction}
It is notable fact that quadratic multi-time Lagrangians are present in most physical domains.
Illustrative examples are present in the theory of elasticity \cite{Olver}, the dynamics of ideal fluids,
    the magnetohydrodynamics \cite{Got-Isen-Mars}, \cite{Holm-Mars-Ratiu} and in
    the theory of bosonic strings \cite{Gi-Mang-Sard}.
This fact encourages the natural attempt of geometrization for quadratic multi-time Lagrangians.
This framework implies, as can be seen below, the introduction of a corresponding
    Riemann-Lagrange geometry on 1-jet spaces.

\section{The generalized multi-time Lagrange space of a quadratic Lagrangian}

Let $U_{(i)}^{(\alpha )}(t^{\gamma },x^{k})$ be a $d$-tensor ({\em distinguished} tensor, in brief)
    on the 1-jet space $J^{1}(T,M)$, and let $F:T\times M\rightarrow \mathbb{R}$ be a smooth function.
We further consider the quadratic multi-time Lagrangians $L:J^{1}(T,M)\rightarrow\mathbb{R}$,
    of the form
\begin{equation}\label{quadratic}L=G_{(i)(j)}^{(\alpha )(\beta )}(t^{\gamma },x^{k})x_{\alpha }^{i}x_{\beta}^{j}+U_{(i)}^{(\alpha )}(t^{\gamma },x^{k})x_{\alpha }^{i}+F(t^{\gamma},x^{k}),\end{equation}
whose fundamental vertical metrical $d$-tensor
$$G_{(i)(j)}^{(\alpha )(\beta )}(t^{\gamma },x^{k})=\frac{1}{2}\frac{\partial^2 L}{\partial x_{\alpha }^{i}\partial x_{\beta }^{j}}$$
is symmetric, of rank $n=\dim M$ and has a constant signature with respect to the indices $i$ and $j$.
By using a semi-Riemannian metric $h=(h_{\alpha\beta }(t^{\gamma }))_{\alpha ,\beta =\overline{1,p}}$ on $T$ and by considering the canonical Kronecker $h$-regular vertical metrical $d$-tensor attached to the Lagrangian function (\ref{quadratic}), given by
$$\mathcal{G}_{(i)(j)}^{(\alpha )(\beta )}(t^{\gamma },x^{k})=\frac{1}{p}h^{\alpha \beta }(t^{\gamma })h_{\mu \nu }(t^{\gamma })G_{(i)(j)}^{(\mu)(\nu )}(t^{\gamma },x^{k}),$$
where $p=\dim T$, we then consider the pair
$$GL(J)=\left( J^{1}(T,M),\mathcal{G}_{(i)(j)}^{(\alpha)(\beta )}(t^{\gamma },x^{k})=h^{\alpha \beta }(t^{\gamma })g_{ij}(t^{\gamma
},x^{k}),\right)$$
which is a generalized multi-time Lagrange space, whose spatial metrical $d$-tensor is given by the formula
\begin{equation}\label{spatial metric}g_{ij}(t^{\gamma },x^{k})=\frac{1}{p}h_{\mu \nu }(t^{\gamma})G_{(i)(j)}^{(\mu )(\nu )}(t^{\gamma },x^{k}).\end{equation}
\begin{defi}
We call the space $GL(J)$ {\em the canonical generalized multi-time Lagrange space associated with the quadratic
Lagrangian function} given by \eqref{quadratic}.
\end{defi}
In order to construct the main Riemann-Lagrange geometric objects of the space $GL(J)$, i.e., its $d$-linear connection, torsions and
curvatures, one needs a nonlinear connection $\Gamma =\left( M_{(\alpha)\beta }^{(i)},N_{(\alpha )j}^{(i)}\right) $ on $J^{1}(T,M)$.
The fundamental vertical metrical $d$-tensor $\mathcal{G}_{(i)(j)}^{(\alpha)(\beta )}(t^{\gamma },x^{k})=h^{\alpha \beta }(t^{\gamma })g_{ij}(t^{\gamma
},x^{k})$, where $g_{ij}$ is given by (\ref{spatial metric}), produces the following natural nonlinear connection \cite[p. 88]{Matrix-Rom}:
$$\begin{array}{cc}M_{(\alpha )\beta }^{(i)}=-H_{\alpha \beta }^{\gamma }x_{\gamma }^{i}, &N_{(\alpha )j}^{(i)}=\Gamma _{jm}^{i}x_{\alpha }^{m}+\dfrac{g^{im}}{2}\dfrac{\partial g_{jm}}{\partial t^{\alpha }},\end{array}$$
where $H_{\alpha \beta }^{\gamma }$ are the Christoffel symbols of the temporal semi-Riemannian metric $h_{\alpha \beta }$,
    \footnote{Throughout the rest of the paper, the constructed geometrical objects will be expressed in local adapted components with respect to previous adapted bases.}
$$\Gamma _{jk}^{i}(t^{\mu },x^{m})=\dfrac{g^{ir}}{2}\left( \dfrac{\partial g_{jr}}{\partial x^{k}}+\dfrac{\partial g_{kr}}{\partial x^{j}}-\dfrac{\partial g_{jk}}{\partial x^{r}}\right)$$
are the generalized Christoffel symbols of the spatial metric $g_{ij}(t^{\gamma },x^{k})$. Let
$$\left\{ \frac{\delta }{\delta t^{\alpha }},\frac{\delta }{\delta x^{i}},\frac{\partial }{\partial x_{\alpha }^{i}}\right\} \subset \mathcal{X}\left(
J^{1}(T,M)\right) \text{ and }\left\{ dt^{\alpha },dx^{i},\delta x_{\alpha}^{i}\right\} \subset \mathcal{X}^{\ast }\left( J^{1}(T,M)\right)$$
be the adapted bases of the nonlinear connection $\Gamma $, where
$$\frac{\delta }{\delta t^{\alpha }}=\frac{\partial }{\partial t^{\alpha }}-M_{(\beta )\alpha }^{(j)}\frac{\partial }{\partial x_{\beta }^{j}},\quad
\frac{\delta }{\delta x^{i}}=\frac{\partial }{\partial x^{i}}-N_{(\beta)i}^{(j)}\frac{\partial }{\partial x_{\beta }^{j}},$$
$$\delta x_{\alpha }^{i}=dx_{\alpha }^{i}+M_{(\alpha )\beta }^{(i)}dt^{\beta}+N_{(\alpha )j}^{(i)}dx^{j}.$$
\section{The Riemann-Lagrange geometry of the space $GL(J)$}
We further adopt the formalism introduced and developed in core seminal research studies (\cite{Asanov},\cite{Matrix-Rom},\cite{Mir-An-Kirk}), and accordingly provide the main results of the Riemann-Lagrange geometry of the generalized multi-time Lagrange space $GL(J)$.
\begin{theo}[the Cartan linear connection]
The canonical Cartan linear connection of the space $GL(J)$ is given by its adapted components
$$C\Gamma =\left( H_{\alpha \beta }^{\gamma },G_{j\gamma}^{k},L_{jk}^{i}=\Gamma _{jk}^{i},C_{j(k)}^{i(\gamma )}=0\right) ,$$
where $G_{j\gamma }^{k}=\dfrac{g^{km}}{2}\dfrac{\partial g_{mj}}{\partial t^{\gamma }}$.
\end{theo}
\proof The formulas which describe the adapted coefficients of the Cartan canonical connection are given by
$$G_{j\gamma }^{k}=\dfrac{g^{km}}{2}\dfrac{\delta g_{mj}}{\delta t^{\gamma }}=\dfrac{g^{km}}{2}\dfrac{\partial g_{mj}}{\partial t^{\gamma }},\quad
L_{jk}^{i}=\dfrac{g^{im}}{2}\left( \dfrac{\delta g_{jm}}{\delta x^{k}}+\dfrac{\delta g_{km}}{\delta x^{j}}-\dfrac{\delta g_{jk}}{\delta x^{m}}\right) =\Gamma _{jk}^{i},$$
\centerline{$C_{j(k)}^{i(\gamma )}=\dfrac{g^{im}}{2}\left( \dfrac{\partial g_{jm}}{\partial x_{\gamma }^{k}}+\dfrac{\partial g_{km}}{\partial x_{\gamma }^{j}}-\dfrac{\partial g_{jk}}{\partial x_{\gamma }^{m}}\right) =0.$ \hfill$\Box$}
\begin{rem}
The generalized Cartan canonical connection $C\Gamma $ of the space $GL(J)$ satisfies the metricity conditions
$$h_{\alpha \beta /\gamma }=h_{\alpha \beta |k}=h_{\alpha \beta}|_{(k)}^{(\gamma )}=0,\qquad g_{ij/\gamma }=g_{ij|k}=g_{ij}|_{(k)}^{(\gamma)}=0,$$
where " $_{/\gamma }$", " $_{|k}$" and " $|_{(k)}^{(\gamma )}$" are the local covariant derivatives produced by the Cartan connection $C\Gamma $.
\end{rem}
\begin{theo}
The generalized multi-time Lagrange space $GL(J)$ is characterized by the following adapted torsion $d$-tensors:
$$T_{\alpha j}^{m}=-G_{j\alpha }^{m},\quad P_{i(j)}^{m(\beta)}=C_{i(j)}^{m(\beta )}=0,\quad P_{(\mu )i(j)}^{(m)\text{ }(\beta )}=\frac{\partial N_{(\mu )i}^{(m)}}{\partial x_{\beta }^{j}}-\delta _{\mu }^{\beta}L_{ij}^{m}=0,$$
$$P_{(\mu )\alpha (j)}^{(m)\text{ }(\beta )}=\frac{\partial M_{(\mu )\alpha}^{(m)}}{\partial x_{\beta }^{j}}-\delta _{\mu }^{\beta }G_{j\alpha
}^{m}+\delta _{j}^{m}H_{\mu \alpha }^{\beta }=-\delta _{\mu }^{\beta}G_{j\alpha }^{m},\quad$$
$$R_{(\mu )\alpha \beta }^{(m)}=\frac{\delta M_{(\mu )\alpha }^{(m)}}{\delta t^{\beta }}-\frac{\delta M_{(\mu )\beta }^{(m)}}{\delta t^{\alpha }},\quad
R_{(\mu )\alpha j}^{(m)}=\frac{\delta M_{(\mu )\alpha }^{(m)}}{\delta x^{j}}-\frac{\delta N_{(\mu )j}^{(m)}}{\delta t^{\alpha }},$$
$$R_{(\mu )ij}^{(m)}=\frac{\delta N_{(\mu )i}^{(m)}}{\delta x^{j}}-\frac{\delta N_{(\mu )j}^{(m)}}{\delta x^{i}},\quad S_{(\mu )(i)(j)}^{(m)(\alpha)(\beta )}=\delta _{\mu }^{\alpha }C_{i(j)}^{m(\beta )}-\delta _{\mu}^{\beta }C_{j(i)}^{m(\alpha )}=0.$$
\end{theo}
Using the general formulas from \cite{Matrix-Rom}, which provide the curvature $d$-tensors of a generalized multi-time Lagrange space, we obtain
the following result:
\begin{theo}
The generalized multi-time Lagrange space $GL(J)$ is characterized by the following adapted curvature d-tensors:
$$H_{\eta \beta \gamma }^{\alpha }=\frac{\partial H_{\eta \beta }^{\alpha }}{\partial t^{\gamma }}-\frac{\partial H_{\eta \gamma }^{\alpha }}{\partial
t^{\beta }}+H_{\eta \beta }^{\mu }H_{\mu \gamma }^{\alpha }-H_{\eta \gamma}^{\mu }H_{\mu \beta }^{\alpha },$$
$$R_{i\beta \gamma }^{l}=\frac{\partial G_{i\beta }^{l}}{\partial t^{\gamma }}-\frac{\partial G_{i\gamma }^{l}}{\partial t^{\beta }}+G_{i\beta}^{m}G_{m\gamma }^{l}-G_{i\gamma }^{m}G_{m\beta }^{l},$$
$$R_{i\beta k}^{l}=\frac{\partial G_{i\beta }^{l}}{\partial x^{k}}-\frac{\partial \Gamma _{ik}^{l}}{\partial t^{\beta }}+G_{i\beta }^{m}\Gamma
_{mk}^{l}-\Gamma _{ik}^{m}G_{m\beta }^{l},$$
$$R_{ijk}^{l}=\frac{\partial \Gamma _{ij}^{l}}{\partial x^{k}}-\frac{\partial\Gamma _{ik}^{l}}{\partial x^{j}}+\Gamma _{ij}^{m}\Gamma _{mk}^{l}-\Gamma
_{ik}^{m}\Gamma _{mj}^{l},$$
$$P_{i\beta (k)}^{l\text{ \ }(\gamma )}=0,\quad P_{ij(k)}^{l\text{ }(\gamma)}=0,\quad S_{i(j)(k)}^{l(\beta )(\gamma )}=0.$$
\end{theo}
\section{Generalized multi-time field theories on the space $GL(J)$}
\subsection{Multi-time gravitational field}

The fundamental vertical metrical d-tensor of the space $GL(J)$ naturally induces the multi-time gravitational $h$-potential $G$, defined %
$$G=h_{\alpha \beta }dt^{\alpha }\otimes dt^{\beta }+g_{ij}dx^{i}\otimes dx^{j}+h^{\alpha \beta }g_{ij}\delta x_{\alpha }^{i}\otimes \delta x_{\beta}^{j}.$$
We postulate that the generalized Einstein equations corresponding to the multi-time gravitational $h$-potential of the space $GL(J)$, are of the form
\begin{equation}\label{Einstein equations}Ric(C\Gamma )-\frac{Sc(C\Gamma )}{2}G=\mathcal{KT},\end{equation}
where $Ric(C\Gamma )$ represents the Ricci $d$-tensor associated with the generalized Cartan connection, $Sc(C\Gamma )$ is the scalar curvature,
    $\mathcal{K}$ is the Einstein curvature scalar and $\mathcal{T}$ is the stress-energy $d$-tensor of matter.

Using now the general formulas from \cite{Matrix-Rom}, we infer the following:
\begin{prop}
The Ricci $d$-tensor $Ric(C\Gamma )$ of the space $GL(J)$ has the following adapted components:
$$\begin{array}{l}\medskip R_{\alpha \beta }:=H_{\alpha \beta }=H_{\alpha \beta \mu }^{\mu},\quad R_{i(j)}^{\;(\alpha )}:=P_{i(j)}^{\;(\alpha )}=-P_{im(j)}^{m\text{ }(\alpha )}=0, \\
    R_{(i)j}^{(\alpha )}:=P_{(i)j}^{(\alpha )}=P_{ij(m)}^{m\text{ }(\alpha)}=0,\quad R_{(i)\beta }^{(\alpha )}:=P_{(i)\beta }^{(\alpha )}=P_{i\beta
(m)}^{m\text{ }(\alpha )}=0,\medskip \\
    R_{(i)(j)}^{(\alpha )(\beta )}:=S_{(i)(j)}^{(\alpha )(\beta)}=S_{i(j)(m)}^{m(\beta )(\alpha )}=0,\quad R_{i\alpha }=R_{i\alpha
m}^{m},\quad R_{ij}=R_{ijm}^{m}.\end{array}$$
\end{prop}

\begin{cor}
The scalar curvature $Sc(C\Gamma )$ of the space $GL(J)$ is given by
$$Sc(C\Gamma )=H+R,$$
where $H=h^{\alpha \beta }H_{\alpha \beta }$ and $R=g^{ij}R_{ij}$.
\end{cor}
Then we can state the main result of the generalized Riemann-Lagrange geometry of the multi-time gravitational field:
\begin{theo}
The global generalized Einstein equations (\ref{Einstein equations}) of the space $GL(J)$, have the local form
$$\left\{\begin{array}{l}\medskip H{_{\alpha \beta }-{\dfrac{H+R}{2}}h_{\alpha \beta }=\mathcal{K}}\mathcal{T}{_{\alpha \beta }} \\
\medskip R{_{ij}-{\dfrac{H+R}{2}}g_{ij}=\mathcal{K}}\mathcal{T}{_{ij}} \\
\medskip {-{\dfrac{H+R}{2}}h^{\alpha \beta }g_{ij}=\mathcal{K}}\mathcal{T}{_{(i)(j)}^{(\alpha )(\beta )}} \\
\begin{array}{lll}\medskip 0=\mathcal{T}_{\alpha i}, & R_{i\alpha }=\mathcal{KT}_{i\alpha }, &0=\mathcal{T}_{(i)\beta }^{(\alpha )} \\
0=\mathcal{T}_{\alpha (i)}^{\text{ \ }\!(\beta )}, & 0=\mathcal{T}_{i(j)}^{\;(\alpha )}, & 0=\mathcal{T}_{(i)j}^{(\alpha )},\end{array}\end{array}\right.$$
where $\mathcal{T}{_{AB}}$, $A,B\in \left\{ \alpha ,\text{ }i,\text{ }_{(i)}^{(\alpha )}\right\}$ are the adapted components of the stress-energy $d$-tensor $\mathcal{T}{.}$
\end{theo}
\subsection{The multi-time electromagnetism}
The multi-time electromagnetic theory of the space $GL(J)$ relies on the \textit{metrical deflection $d$-tensors}
$$\begin{array}{c}\medskip {D_{(i)j}^{(\alpha )}=}\left[ \mathcal{G}_{(i)(m)}^{(\alpha )(\mu)}x_{\mu }^{m}\right] _{|j}{=-{\dfrac{h^{\alpha \mu }}{2}}{\dfrac{\partial g_{ij}}{\partial t^{\mu }}},} \\
d_{(i)(j)}^{(\alpha )(\beta )}=\left[ \mathcal{G}_{(i)(m)}^{(\alpha )(\mu)}x_{\mu }^{m}\right] |_{(j)}^{(\beta )}{=}h^{\alpha \beta }g_{ij}.\end{array}$$
This yields the {\em electromagnetic $2$-form} of the space $GL(J)$, via:
$$F=F_{(i)j}^{(\alpha )}\delta x_{\alpha }^{i}\wedge dx^{j}+f_{(i)(j)}^{(\alpha )(\beta )}\delta x_{\alpha }^{i}\wedge \delta x_{\beta }^{j},$$
where
$$F_{(i)j}^{(\alpha )}={{\dfrac{1}{2}}\left[ D_{(i)j}^{(\alpha)}-D_{(j)i}^{(\alpha )}\right] =0},\qquad f_{(i)(j)}^{(\alpha )(\beta )}={{\dfrac{1}{2}}\left[ d_{(i)(j)}^{(\alpha )(\beta )}-d_{(j)(i)}^{(\alpha)(\beta )}\right] =0}.$$
Since $F=0$, we infer that the multi-time electromagnetic theory of the space $GL(J)$ is formally trivial.\medskip

\noindent\textbf{Acknowledgements.} The author thanks to Professor Vladimir Balan, whose advice helped to improve this paper.
\noindent Mircea Neagu\\
Transilvania University of Bra\c{s}ov,\\
Department of Mathematics and Informatics,\\
Blvd. Iuliu Maniu 50, 500091 Bra\c{s}ov, Romania.\\
E-mail: mircea.neagu@unitbv.ro\\\\

\begin{thebibliography}{9}

\bibitem{Asanov}G.S. Asanov, {\em Jet extension of Finslerian gauge approach}, Fortschritte der Physik, vol. {\bf 38}, 8 (1990), 571-610.

\bibitem{Gi-Mang-Sard}G. Giachetta, L. Mangiarotti, G. Sardanashvily, {\em Covariant Hamitonian field theory}, arXiv:hep-th/9904062, 1999.

\bibitem{Got-Isen-Mars}M. Gotay, J. Isenberg, J.E. Marsden, R. Montgomery, {\em Momentum maps and classical relativistic fields.
    Part I: Covariant field theory}, arXiv:physics/9801019v2 [math-ph], 2004.

\bibitem{Holm-Mars-Ratiu}D.D. Holm, J.E. Marsden, T.S. Ra\c{t}iu, {\em The Euler-Poincar\'{e} equations and semidirect products
    with applications to continuum theories, }arXiv:chao-dyn/9801015, 1998.

\bibitem{Mir-An-Kirk}R. Miron, M.S. Kirkovits, M. Anastasiei, {\em A geometrical model for variational problems of multiple integrals},
    Proc. of Conf. on Diff. Geom. and Appl., Dubrovnik, Yugoslavia, June 26-July {\bf 3} (1988), 8-25.

\bibitem{Matrix-Rom}M. Neagu, {\em Riemann-Lagrange Geometry on 1-Jet Spaces}, Matrix Rom, Bucharest, 2005.

\bibitem{Olver}P.J. Olver, {\em Canonical elastic moduli}, Journal of Elasticity {\bf 19} (1988), 189-212.

\end{thebibliography}
\end{document}